\documentclass{amsproc}
\usepackage{graphicx}
\usepackage{caption}
\usepackage{subfigure}

\usepackage{amssymb}
\usepackage{epstopdf}
\usepackage{enumitem}
\usepackage{amsthm,amsmath}
\usepackage{amsrefs}
\usepackage{framed}
\usepackage{hyperref}
\usepackage{xcolor}

\hypersetup{
    colorlinks,
    linkcolor={red!70!black},
    citecolor={blue!70!black},
    urlcolor={blue!50!black}
}
\usepackage{tikz}
\usetikzlibrary{shapes,arrows,calc,positioning,scopes}
\usepackage{verbatim}

\usepackage{multicol}

\theoremstyle{plain}

\newtheorem*{qst}{Qualitative Stability Theorem}
\theoremstyle{definition}

\newtheorem*{defcyc}{Definition of $k$-cycle}
\newtheorem*{defcycs}{Definition of simple $k$-cycle}

\newtheorem{exmp}{Example}
\theoremstyle{remark}

\newcommand{\bb}[0]{\begin{bmatrix}}
\newcommand{\eb}[0]{\end{bmatrix}}

\renewcommand{\Re}[0]{\mathbb R}

\newcommand{\V}[0]{\mathcal V}
\newcommand{\E}[0]{\mathcal E}

\newcommand{\be}[0]{\begin{equation}}
\newcommand{\ee}[0]{\end{equation}}
\newcommand{\ben}[0]{\begin{equation*}}
\newcommand{\een}[0]{\end{equation*}}

\title{Sign Stability via Root Locus Analysis}
\author{Travis E. Gibson}
\address{Research Fellow in Medicine, Harvard Medical School and Brigham and Women's Hospital}
\email{tgibson@mit.edu}
\date{\today}                                  
\begin{document}
\maketitle

\noindent{\bf Abstract--}With the rise of network science old topics in ecology and economics are resurfacing.  One such topic is structural stability (often referred to as qualitative stability or sign stability). A system is deemed structurally stable if the system remains stable for all possible parameter variations so long as the parameters do not change sign. This type of stability analysis is appealing when studying real systems as the underlying stability result only requires the scientist or engineer to know the sign of the parameters in the model and not the specific values. The necessary and sufficient conditions for qualitative stability however are opaque. In order to shed light on those conditions root locus analysis is employed. This technique allows us to illustrate the necessary conditions for qualitative stability.

\section*{Introduction}
Knowing when the eigenvalues of a real matrix are located in the closed left half of the complex plane is of fundamental importance when studying the stability of linear systems. When the above criterion is met for $A$ in the linear system $\dot x = Ax$, all trajectories are bounded above by an exponentially decreasing function. Often in real systems that are known to follow a linear dynamical process, the exact values of the entries in the  $A$ matrix may not be known. Most likely some information will be known however. For instance, what if only the sign of the entries in the $A$ matrix were known, then can one still say something about the stability of said system? This is precisely the question that James Quirk and Richard Ruppert sought to answer when they introduced the notion of qualitative stability \cite{quirk1965qualitative}.  A matrix was deemed qualitatively stable if all matrices with the same sign pattern had eigenvalues in the closed left half of the complex plane. This stability condition is also referred to as sign stability or structural stability. 

While first introduced in the field of economics, qualitative stability was quickly introduced to ecologist through the works of Clark Jeffries \cite{Jeffries:1974aa} and Robert May \cite{may1973qualitative}. The field of complex systems has taken off recently with renewed interest in graph theory \cite{bollobas1998modern}, its applications in dynamical systems \cite{fax_tac_04,olfati2004consensus,dor13}, and the newly established field of network science \cite{new03,alb02review}. It is because of this renewed interest that we wish to give further insights into the necessary and sufficient conditions for the sign stability of a real matrix.\footnote{We note that the use of graphs to illustrate and understand dynamic processes over a network is an almost 60 year old idea \cite{coa_ire_59,mas_ire_53,mas_ire_56}} In this short note we will use the tool of root-locus, first developed in the field of electrical engineering \cite{5059708,5060121}, so as to give motivation for the necessary and sufficient conditions for sign stability. The main insight comes by using the fact that all $k$-cycles of length 3 or greater will always have an asymptote of the root locus going into the right half of the complex plane.

\section*{Sign Stability}
We now formally define our dynamics of interest. Consider a dynamical system with the state vector $x= [x_1, x_2,\ldots, x_n]^T: \Re \to \Re^n $ defined by the linear differential equation
\begin{equation}\label{eq:linA}
\dot x(t) = A x(t)
\end{equation}
where $A\in \Re^{n\times n}$ whose components are indexed as $a_{ij}=[A]_{ij}$, and $t\in\Re$ is time. 
The following is a slightly weakened version of the original result {\cite[{Theorem 5}]{quirk1965qualitative}} which was proven to only hold generically \cite{yamada1987generic}.

\begin{qst}\label{thm:signstab} The matrix $A\in\Re^{n\times n}$ is Hurwitz for {\em almost} all values of $a_{ij}$ with constant sign if and only if 
\begin{enumerate}[label={\normalfont (\roman*)}]
\item $a_{ii} \leq 0 $ for all $i$.
\item $a_{ii} \neq 0$ for at least one $i$.
\item $a_{ij}a_{ji}\leq 0$, $i\neq j$.
\item No $k$-cycles of length $k\geq 3$.
\item $\text{det}({A})\neq 0$.
\end{enumerate}
\end{qst}
\noindent Now consider two matrices $B$ and $C$ with sign patterns as follows
\ben
\textrm{sign}({B})=\bb - & + & + & + \\ 
- & - & 0 & 0 \\
- & 0 & - & 0 \\
- & 0 & 0 &  - \eb \quad \text{and} \quad \textrm{sign}({C})=\bb - & + & + & + \\ 
- & - & + & + \\
- & - & - & + \\
- & - & - &  - \eb.
\een
Using the above result, we know that any matrix with the same sign pattern as $B$ will almost always be stable, while there are many matrices with the same sign pattern as $C$ that have at least one eigenvalue in the the right half of the complex plane \cite[Equations (9) and (10)]{may1973qualitative}. Before discussing the conditions for sign stability in the theorem there are two preliminary results that are needed. The first is the equivalence between a directed graph with a weighted adjacency matrix and a network of linear operators. The second preliminary discussion involves the rules for placing the asymptotes of the root locus for our linear operators of interest.

\section*{Equivalence of Directed Graphs and Linear Networks}
The dynamics in \eqref{eq:linA} have equivalent representations as both a digraph and a network of linear operators. The latter is discussed first. The dynamics in \eqref{eq:linA} can be represented as $n$ separate dynamical systems
\ben
\dot x_i(t) = a_{ii} x_i + \sum_{j\neq i} a_{ij} x_j(t), \quad i\in \{1,2,\dots,n\}.
\een
Defining $s\triangleq \frac{d}{dt}$ as the differential operator, the dynamics have an equivalent representation as 
\be \label{eq:s}
x_i(t) = \frac{1}{s-a_{ii}} \sum_{j\neq i} a_{ij} x_j(t), \quad i\in \{1,2,\ldots,n\}.
\ee
The dynamics in \eqref{eq:s} can be represented as a network where each node is a linear operator of the form  $\frac{1}{s-a_{ii}}$ and the corresponding incoming and outgoing edges have linear multipliers of the form $a_{ij}$ and $a_{ji}$ where $j\neq i$. All edges out of node $i$ cary the flow of the state variable $x_i(t)$. The dynamics in \eqref{eq:linA}  also have an equivalent digraph $\mathcal D(\mathcal V,\mathcal E,\mathcal A)$ representation, where $\mathcal V = \{1,2,\ldots,n\}$ is the vertex set, the directed edges are defined by the ordered pairs ${(i, j)\in \E \subset{\V}\times{\V}}$ and $\mathcal A^T=A$ is the weighted adjacency matrix. An element $(i,j)\in \E$ if and only if there is a directed edge from vertex $i$ to vertex $j$, i.e. $\alpha_{ij}=[\mathcal A]_{ij}$ is non zero. The two equivalent representation are now compared for a graph and network containing a 3-cycle, but first we give the formal definition of a $k$-cycle in terms of a digraph.

\begin{defcyc}
A $k$-cycle on a digraph is a sequences of $k+1$ vertices beginning and ending with the same vertex with no other vertices repeated, with a directed edge existing to complete a walk of the sequence.
\end{defcyc}

\begin{exmp}\label{ex:1}
Consider a third order system satisfying \eqref{eq:linA} where
\be
A= \bb 
	a_{11} & 0 & a_{13} \\ 
       a_{21} & a_{22} & 0\\
       0 & a_{32} &a_{33}
      \eb.
\ee
A digraph and linear network for this example are shown in Figure \ref{fig:3cyc}. The first thing to notice is that the self loops of the digraph are absorbed into the nodes of the network  in the form of a linear operator. The second noticeable difference in these two representations is the transposition of the indexing that comes via our definition $\mathcal A=A^T$, and as a result $\alpha_{ij}=a_{ji}$. 
\end{exmp}

\tikzstyle{block} = [draw, rectangle, minimum height=2em]
\tikzstyle{circ} = [draw, circle]
\tikzstyle{ellip} = [draw, ellipse]
\tikzstyle{li} = [draw, line]
\tikzstyle{input} = [coordinate]
\tikzstyle{output} = [coordinate]
\tikzstyle{pinstyle} = [pin edge={to-,thin,black}]

\begin{figure}[t]
\centering
\subfigure[Digraph with 3-cycle] {
\begin{tikzpicture}[->,>=stealth',shorten >=1pt,auto,node distance=1.3cm,main node/.style={circle,draw},font=\footnotesize]
  \node[main node] (1) {1};
  \node[main node] (2) [right of=1] {2};
  \node[main node] (3) [right of=2] {3};
  \path
    (1) edge [loop above] node {$\alpha_{11}$} (1)
         edge [right] node [above]{$\alpha_{12}$} (2)
    (2) edge [right] node [above] {$\alpha_{23}$} (3)
         edge [loop above] node {$\alpha_{22}$} (2)
    (3) edge [bend left] node [below ]{$\alpha_{31}$} (1)
         edge [loop above] node {$\alpha_{33}$} (3);
\end{tikzpicture}
}
\subfigure[Linear network of 3-cycle]{
\begin{tikzpicture}[->,>=stealth',shorten >=1pt,auto, node distance=1.5cm,font=\footnotesize]

    \node [circ, name=a11] {$\frac{1}{s-a_{11}}$};
    \node [block, right of=a11] (a12){$a_{21}$};
    \node [circ, right of=a12] (a22) {$\frac{1}{s-a_{22}}$};
    \node [block, right of=a22] (a23){$a_{32}$};
    \node [circ, right of=a23] (a33) {$\frac{1}{s-a_{33}}$};
    \node [block, below of=a22] (a31){$a_{13}$};
     \draw [->] (a11) -- node[name=x1] {$x_1$} (a12);
     \draw [->] (a12) -- node{}(a22);
     \draw [->] (a22) -- node[name=x2] {$x_2$} (a23);
     \draw [->] (a23) -- node{}(a33);
     \draw [->] (a33) |- node[above left] {$x_3$} (a31);
     \draw [->] (a31) -| node{}(a11);
     
\end{tikzpicture} }
\caption{Visualization of the dynamics from Example \ref{ex:1} in (a) graph and (b) network of linear operators. Recall that the adjacency matrix is defined as $\mathcal A=A^T$ and the linear transfer function network is defined as $sx=Ax$.  Also, a reminder on or  notation for indexing matrices, $a_{ij}=[A]_{ij}$ and $\alpha_{ij}=[\mathcal A]_{ij}$.}
\label{fig:3cyc}
\end{figure}
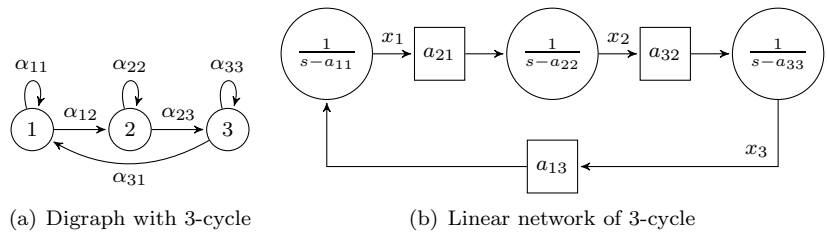

\section*{Simple $k$-cycles}
\begin{figure*}[B!]
\centering
\subfigure[$k$-cycle]{
\begin{tikzpicture}[->,>=stealth',shorten >=1pt,auto, node distance=1.8cm,font=\footnotesize]
   \node [circ, name=a11] {$\frac{1}{s-a_{11}}$};
    \node [block, right of=a11] (a12){$a_{21}$};
    \node [circ, right of=a12] (a22) {$\frac{1}{s-a_{22}}$};
    \node [block, right of=a22] (a23){$a_{32}$};
     \node [right of=a23] (add){$\ldots$};
    \node [block, right of=add] (an1n){$a_{k(k-1)}$};
    \node [circ, right of=an1n] (ann) {$\frac{1}{s-a_{kk}}$};
    \node [block, below of=a23] (an1){$a_{1k}$};
    \draw [->] (a11) -- node{$x_1$} (a12);
    \draw [->] (a12) -- node{}(a22);
    \draw [->] (a22) -- node{$x_2$} (a23);
     \draw [-] (a23) -- (add);
     \draw [->] (add) -- node{$x_{k-1}$} (an1n);
     \draw [->] (an1n) -- (ann);
     \draw [->] (ann) |- node [above left]{$x_k$} (an1);
     \draw [->] (an1) -| node{}(a11);
\end{tikzpicture}}

\subfigure[$k$-cycle with one feedback gain]{
\begin{tikzpicture}[->,>=stealth',shorten >=1pt,auto, node distance=1.8cm,font=\footnotesize]
    \node [circ, name=a11] {$\frac{1}{s-a_{11}}$};
    \node [circ, right of=a11] (a22) {$\frac{1}{s-a_{22}}$};
     \node [right of=a22] (add){$\ldots$};
    \node [circ, right of=add] (ann) {$\frac{1}{s-a_{kk}}$};
     \draw [->] (a11) -- node{$y_1$}(a22);
     \draw [-] (a22) -- node[name=bl]{$y_2$}(add);
     \node [block, below of = bl] (an1){$a_{21}a_{32}\cdots a_{k(k-1)}a_{1k}$};
     \draw [->] (add) --  node{$y_{k-1}$}(ann);
     \draw [->] (ann) |-  node[above left]{$y_k$}(an1);
     \draw [->] (an1) -| node{}(a11);
 \end{tikzpicture}}

\subfigure[Compact representation of $k$ cycle with one feedback gain $\ell$]{
\begin{tikzpicture}[->,>=stealth',shorten >=1pt,auto, node distance=1.8cm,font=\footnotesize]
   \node [ellip, name=a11] {$\frac{1}{(s-a_{11})(s-a_{22})\cdots(s-a_{kk})}$};
   \node [block, below of = a11] (an1){$\ell$};
   \draw [->] (a11) -| (3,-1) |- node[above left]{$y_k$}(an1);
      \draw [->] (an1) -| (-3,-1) |- (a11);
\end{tikzpicture} }
\caption{Visualization of the state space transformations from \eqref{eq:kx} to \eqref{eq:kk}.}
\label{fig:comp}
\end{figure*}
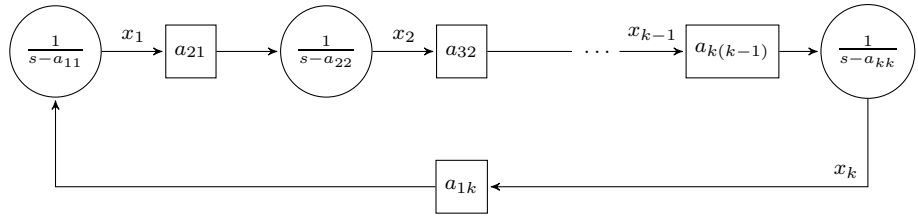
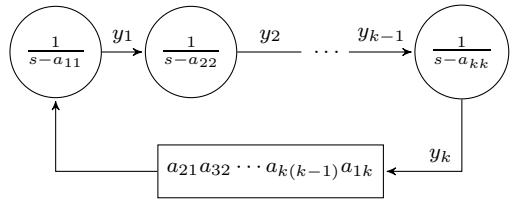
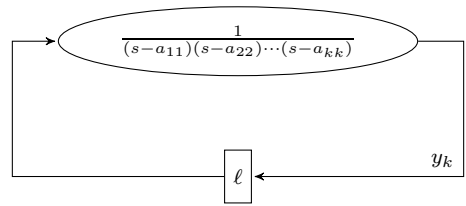

In this section we are only interested in linear dynamics of the form \eqref{eq:linA} for when $A\in \Re^{k\times k}$ and the associated digraph contains only a $k$-cycle and 1-cycles. This class of digraphs is now formally defined.
\begin{defcycs}
The digraph $\mathcal D(\mathcal V,\mathcal E,\mathcal A)$ with $\mathcal V=\{1,2,\ldots,k\}$ is a simple $k$-cycle if $\mathcal D$ contains a single $k$-cycle and all other edges in $\mathcal E$ are of the form ($i,i$), $i\in \mathcal V$. 
\end{defcycs}\noindent The indices of a simple $k$-cycle can always be permuted so that the equivalent linear system satisfies \eqref{eq:linA} where
\be\label{eq:Ak}
A= \bb 
	a_{11} & 0 & \cdots & 0 & a_{1k} \\ 
       a_{21} & a_{22} & 0 & \cdots &0\\
       0 & \ddots & \ddots & \ddots & \vdots  \\
              \vdots & \ddots & \ddots & \ddots& 0   \\
       0 & \cdots  & 0 & a_{k(k-1)} & a_{kk}  
      \eb.
\ee
The dynamics of a simple $k$-cycle therefore can always be written in terms of a linear network 
\be\label{eq:kx}
\begin{split}
x_1(t) &=  \frac{1}{s-a_{11}} a_{1k} x_{k}(t) \\
x_i(t)  &=  \frac{1}{s-a_{ii}}  a_{i(i-1)} x_{i-1}(t), \quad i\in\{2,3,\ldots,k\}. 
\end{split}
\ee
Do to linearity of the operators $\frac{1}{s-a_{ii}}$, all of the terms $a_{ij}$, $i\neq j$, can be passed into a single variable  $\ell=a_{1k}\prod_{i=1}^{k-1} a_{{(i+1)}i}$ which we will refer to as the {\em loop gain}. Then a new state variable $y=[y_1,y_2,\ldots y_k]^T$ can defined as $y=Tx$ for some $T\neq 0$ such that $y$ satisfies the following linear relation
\be\label{eq:ky}
\begin{split}
y_1(t) &=  \frac{1}{s-a_{11}} \ell y_{k}(t) \\
y_i(t)  &=  \frac{1}{s-a_{ii}}   y_{i-1}(t), \quad i\in\{2,3,\ldots,k\}. 
\end{split}
\ee
By successive substitution of the relations in \eqref{eq:ky} for $i\neq k$ we have \be\label{eq:kk} y_k(t) =\frac{\ell}{p(s)}y_k(t)\ee  where $p(s)=\prod_{i=1}^{k}{(s-a_{ii})}$. In order to facilitate the transition from \eqref{eq:kx} to \eqref{eq:kk} the networks associated with these relations are shown in Figure \ref{fig:comp}.

In order for \eqref{eq:kk} to be satisfied the $k$ roots of
\be\label{eq:ps}
1-\frac{\ell}{p(r)}=0
\ee
 where ${r\in \mathbb C}$, must coincide with the eigenvalues of the characteristic solution to $y(t)$. By construction, the roots of \eqref{eq:ps} are equivalent to the eigenvalues of $A$ in \eqref{eq:Ak} as well. We are now ready to introduce the tools from root locus analysis. Next we are going to illustrate how the roots of \eqref{eq:ps} change as a function $\ell$.

\section*{Root-Locus for $k$-cycles}
\begin{figure*}[b]
\centering
\subfigure[2-cycle with $\ell<0$]{
\begin{tikzpicture}[->,>=stealth',auto, node distance=1.8cm]
   \draw [->] (0,0) -- (3,0) node[below]{$\Re$};
    \draw [->] (2.5,-1.5) -- (2.5,1.5) node[right]{$\mathbb I$};
    \draw [<->,thick,dashed] (1.5,-1.5) -- (1.5,1.5);
\end{tikzpicture}} \quad
\subfigure[3-cycle with $\ell <0$]{
\begin{tikzpicture}[->,>=stealth',shorten >=1pt,auto, node distance=1.8cm]
   \draw [->] (0,0) -- (3,0) node[below]{$\Re$};
    \draw [->] (2.5,-1.5) -- (2.5,1.5) node[right]{$\mathbb I$};
        \draw [<-,thick,dashed] (0,0) -- (1.5,0);
        \draw [->,thick,dashed] (1.5,0) -- (1.5+0.86,1.5);
        \draw [->,thick,dashed] (1.5,0) -- (1.5+0.86,-1.5);
\end{tikzpicture}} \quad
\subfigure[4-cycle with $\ell<0$]{
\begin{tikzpicture}[->,>=stealth',shorten >=1pt,auto, node distance=1.8cm]
    \draw [->] (-.1,0) -- (3,0) node[below]{$\Re$};
    \draw [->] (2.5,-1.5) -- (2.5,1.5) node[right]{$\mathbb I$};
        \draw [->,thick,dashed] (1.5,0) -- (1.5+1.5,1.5);
        \draw [->,thick,dashed] (1.5,0) -- (1.5+1.5,-1.5);
        \draw [->,thick,dashed] (1.5,0) -- (0,1.5);
        \draw [->,thick,dashed] (1.5,0) -- (0,-1.5);
\end{tikzpicture}}

\subfigure[2-cycle with $\ell>0$]{
\begin{tikzpicture}[->,>=stealth',shorten >=1pt,auto, node distance=1.8cm]
   \draw [->] (0,0) -- (3,0) node[below]{$\Re$};
    \draw [->] (2.5,-1.5) -- (2.5,1.5) node[right]{$\mathbb I$};
    \draw [<->,thick,dashed] (0,0) -- (3,0);
\end{tikzpicture}}
\quad
\subfigure[3-cycle with $\ell>0$]{
\begin{tikzpicture}[->,>=stealth',shorten >=1pt,auto, node distance=1.8cm]
   \draw [->] (0,0) -- (3,0) node[below]{$\Re$};
    \draw [->] (2.5,-1.5) -- (2.5,1.5) node[right]{$\mathbb I$};
        \draw [->,thick,dashed] (1.5,0) -- (3,0);
        \draw [->,thick,dashed] (1.5,0) -- (1.5-0.86,1.5);
        \draw [->,thick,dashed] (1.5,0) -- (1.5-0.86,-1.5);
\end{tikzpicture}}
\quad
\subfigure[4-cycle with $\ell>0$]{
\begin{tikzpicture}[->,>=stealth',shorten >=1pt,auto, node distance=1.8cm]
    \draw [->] (-.1,0) -- (3,0) node[below]{$\Re$};
    \draw [->] (2.5,-1.5) -- (2.5,1.5) node[right]{$\mathbb I$};
           \draw [<->,thick,dashed] (1.5,-1.5) -- (1.5,1.5);
           \draw [<->,thick,dashed] (0,0) -- (3,0);
\end{tikzpicture}}
\caption{Asymptotes for  simple 2-cycle, 3-cycle, and 4-cycle with both negative and positive feedback.}
\label{fig:rl}
\end{figure*}
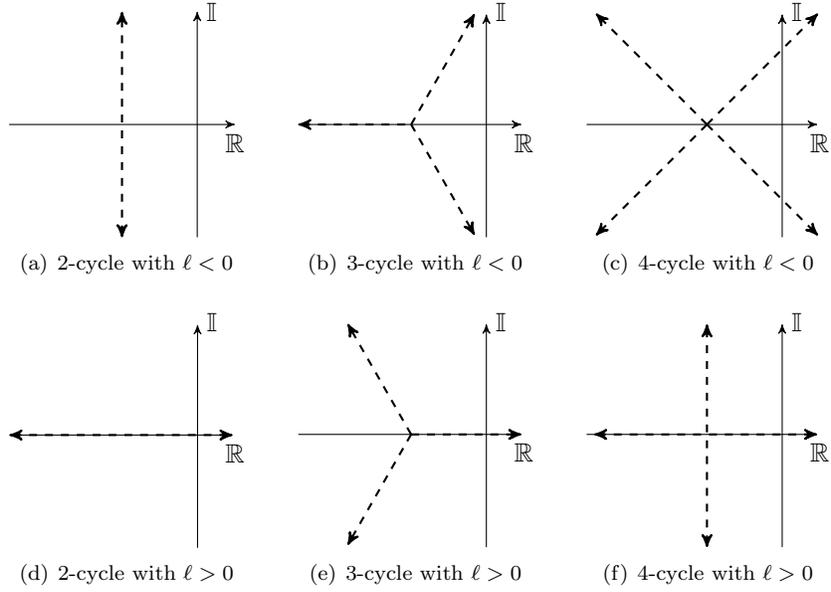

The asymptotes of the root locus for \eqref{eq:ps} as a function of $\ell\to-\infty$ and $\ell\to+\infty$ are now introduced \cite[{Chapter 8}]{nise2010control}. First, given that the expression $\ell/p(s)$ has a zero order numerator, there will be no finite eigenvalues in the limit as  $\ell$ tends to either $-\infty$ or $+\infty$. The intersection of the asymptotes with the real axis in both limiting cases is denoted as
\be\label{eq:sigma}
\sigma =\frac{1}{k}\sum_{i=1}^k a_{ii}.
\ee 
The angle of the $k$ asymptotes in the complex plane, starting from the point $\sigma$, are now defined as
\be\label{eq:thetas}
\begin{split}
\theta_i^-  = & \frac{2i+1}{k}\pi   \\ 
\theta_i^+ = & \frac{2i}{k}\pi
\end{split}
\ee
with $i=0,\pm1,\pm2,\ldots$,  for the cases  $\ell\to-\infty$ and  $\ell\to+\infty$ respectively. The asymptotes for a simple $2$-cycle, $3$-cycle, and $4$-cycle with both negative and positive loop gains are illustrated in Figure \ref{fig:rl}. Three key features of the root locus for simple $k$-cycles are now stated.\footnote{Condition (a) follows by direct inspection of \eqref{eq:sigma}. Conditions (b) follows by observing from \eqref{eq:thetas}  that  $\theta_i^+=0$ for all positive integers $k$ when $i=0$, and  that  $\theta_i^-<\pi/2$ for all integers $k\geq 3$ when $i=0$. Condition (c) follows from the proof of (b) and by the fact that for $k=2$, $\theta^-_i\equiv \pi/2\ (\textrm{mod}\ 2\pi)$ for all even integers $i$ (including 0) and   $\theta^-_i\equiv3\pi/2\ (\textrm{mod}\ 2\pi)$ for all odd integers $i$. }
\begin{enumerate}[label=(\alph*)]
\item If  $a_{ii}\leq0$ for all $i$ and  $a_{jj}<0$ for at least one $j$ then $\sigma<0$.
\item No simple $k$-cycle of length $k\geq 3$ has asymptotes that remain in the closed left half of the complex plane. 
\item For $2$-cycles with $\ell>0$ the asymptotes are at $0$ and $\pi$ and  for $2$-cycles with $\ell<0$ the asymptotes are $\pi/2$ and $3\pi/2$
\end{enumerate}
Thus for simple $k$-cycles of length $k>1$ structural stability is only possible if $k=2$, $\ell\leq0$, and $\sigma<0$.

\section*{Discussion  and Conclusions}
We now discuss how the three properties of simple $k$-cycles illustrate the conditions for qualitative stability. Property (a) shows that with conditions (i) and (ii) of the theorem satisfied $\sigma<0$ and thus the root-loci intersect the real axis in the left-half of the complex plane. Condition (iv) of the theorem is identical to (b) and illustrates the fact that for cycles of length 3 or greater there always exists an $\ell$ such that an eigenvalue is in the right-half plane. Condition (iii) of the theorem, $a_{ij}a_{ji}\leq 0$, is equivalent to $\ell\leq 0$ when $k=2$ and thus is identical to property (c). We hope that this short discussion has illuminated the conditions for sign stability of a matrix. Namely the rather non-intuitive conditions of no $k$-cycles greater than 2, and $a_{ij} a_{ji}\leq 0,\ i \neq j$.

%
%


\bibliography{../master}

\begin{bibdiv}
\begin{biblist}

\bib{alb02review}{article}{
      author={Albert, R{\'e}ka},
      author={Barab{\'a}si, Albert-L{\'a}szl{\'o}},
       title={Statistical mechanics of complex networks},
        date={2002},
     journal={Reviews of modern physics},
      volume={74},
      number={1},
       pages={47},
}

\bib{bollobas1998modern}{book}{
      author={Bollob{\'a}s, B{\'e}la},
       title={Modern graph theory},
   publisher={Springer},
        date={1998},
      volume={184},
}

\bib{coa_ire_59}{article}{
      author={Coates, C.L.},
       title={Flow-graph solutions of linear algebraic equations},
        date={1959Jun},
     journal={Circuit Theory, IRE Transactions on},
      volume={6},
      number={2},
       pages={170\ndash 187},
}

\bib{dor13}{article}{
      author={D{\"o}rfler, Florian},
      author={Chertkov, Michael},
      author={Bullo, Francesco},
       title={Synchronization in complex oscillator networks and smart grids},
        date={2013},
     journal={Proceedings of the National Academy of Sciences},
      volume={110},
      number={6},
       pages={2005\ndash 2010},
}

\bib{5059708}{article}{
      author={Evans, Walter~R.},
       title={Graphical analysis of control systems},
        date={1948},
        ISSN={0096-3860},
     journal={American Institute of Electrical Engineers, Transactions of the},
      volume={67},
      number={1},
       pages={547\ndash 551},
}

\bib{5060121}{article}{
      author={Evans, Walter~R.},
       title={Control system synthesis by root locus method},
        date={1950},
     journal={American Institute of Electrical Engineers, Transactions of the},
      volume={69},
      number={1},
       pages={66\ndash 69},
}

\bib{fax_tac_04}{article}{
      author={Fax, J.A.},
      author={Murray, R.M.},
       title={Information flow and cooperative control of vehicle formations},
        date={2004Sept},
        ISSN={0018-9286},
     journal={Automatic Control, IEEE Transactions on},
      volume={49},
      number={9},
       pages={1465\ndash 1476},
}

\bib{Jeffries:1974aa}{article}{
      author={Jeffries, Clark},
       title={Qualitative stability and digraphs in model ecosystems},
        date={1974},
     journal={Ecology},
      volume={55},
      number={6},
       pages={1415\ndash 1419},
}

\bib{mas_ire_53}{article}{
      author={Mason, S.J.},
       title={Feedback theory-some properties of signal flow graphs},
        date={1953},
        ISSN={0096-8390},
     journal={Proceedings of the IRE},
      volume={41},
      number={9},
       pages={1144\ndash 1156},
}

\bib{mas_ire_56}{article}{
      author={Mason, S.J.},
       title={Feedback theory-further properties of signal flow graphs},
        date={1956},
        ISSN={0096-8390},
     journal={Proceedings of the IRE},
      volume={44},
      number={7},
       pages={920\ndash 926},
}

\bib{may1973qualitative}{article}{
      author={May, Robert~M},
       title={Qualitative stability in model ecosystems},
        date={1973},
     journal={Ecology},
       pages={638\ndash 641},
}

\bib{new03}{article}{
      author={Newman, M.},
       title={The structure and function of complex networks},
        date={2003},
     journal={SIAM Review},
      volume={45},
      number={2},
       pages={167\ndash 256},
}

\bib{nise2010control}{book}{
      author={Nise, N.S.},
       title={Control systems engineering},
   publisher={Wiley},
        date={2010},
}

\bib{olfati2004consensus}{article}{
      author={Olfati-Saber, Reza},
      author={Murray, Richard~M},
       title={Consensus problems in networks of agents with switching topology
  and time-delays},
        date={2004},
     journal={Automatic Control, IEEE Transactions on},
      volume={49},
      number={9},
       pages={1520\ndash 1533},
}

\bib{quirk1965qualitative}{article}{
      author={Quirk, James},
      author={Ruppert, Richard},
       title={Qualitative economics and the stability of equilibrium.},
        date={1965},
     journal={Review of Economic Studies},
      volume={32},
      number={4},
}

\bib{yamada1987generic}{article}{
      author={Yamada, Takeo},
       title={Generic matrix sign-stability.},
        date={1987},
     journal={Can. Math. Bull.},
      volume={30},
      number={3},
       pages={370\ndash 376},
}

\end{biblist}
\end{bibdiv}

\end{document}